\documentclass[preprint]{article}

\usepackage{amssymb}

\newtheorem {definition}{\bf Definition}
\newtheorem{proposition}{\bf Proposition}
\newtheorem {lemma} {\bf Lemma}
\newtheorem{corollary} {\bf Corollary}
\newtheorem {theorem} {\bf Theorem}
\newtheorem{example}{\bf Example}

\title{Calculation of the Number of all Pairs of Disjoint S-permutation Matrices}

\author{Krasimir Yordzhev}
\date{}

\begin{document}
\maketitle

\begin{center}{\it Faculty of Mathematics and Natural Sciences,
 South-West University,
2700 Blagoevgrad, Bulgaria}\\
Email address: yordzhev@swu.bg

\end{center}

\textbf{Abstract:}
\emph{The concept of S-permutation matrix is considered. A general formula for counting all disjoint pairs of $n^2 \times n^2$ S-permutation matrices  as a function of the positive integer $n$ is formulated and proven in this paper. To do that, the graph theory techniques have been used. It has been shown that to count the number of disjoint pairs of $n^2 \times n^2$ S-permutation matrices, it is sufficient to obtain some numerical characteristics of all $n\times n$ bipartite graphs.}



\section{Introduction and notation}
The text of the present paper is essentially extended and modified to the previous articles of the author ~\cite{Yordzhev2013,Yordzhev201447}, while we have corrected all inaccuracies and calculation errors.

Let $n$ be a positive integer. By $[n]$ we denote the set $[n] =\left\{ 1,2,\ldots ,n\right\}$.

A \emph{binary} (or   \emph{boolean}, or (0,1)-\emph{matrix})
 is a matrix all of whose elements belong to the set $\mathfrak{B}
=\{ 0,1 \}$. With $\mathfrak{B}_n$ we will denote the set of all  $n
\times n$  binary matrices.

Let $n$ be a positive integer and let $A\in \mathfrak{B}_{n^2}$ be a $n^2 \times n^2$ binary matrix.  With the help of $n - 1$ horizontal lines and $n - 1$ vertical lines $A$ has been separated into $n^2$ of number non-intersecting $n\times n$ square sub-matrices $A_{kl}$, $1\le k,l\le n$, e.i.
\begin{equation}\label{matrA}
A =
\left[
\begin{array}{cccc}
A_{11} & A_{12} & \cdots & A_{1n} \\
A_{21} & A_{22} & \cdots & A_{2n} \\
\vdots & \vdots & \ddots & \vdots \\
A_{n1} & A_{n2} & \cdots & A_{nn}
\end{array}
\right] .
\end{equation}

The sub-matrices $A_{kl}$, $1\le k,l\le n$ will be called \emph{blocks}.

 A matrix $A\in \mathfrak{B}_{n^2}$ is called an \emph{S-permutation} if in each row, each column, and each block of $A$ there is exactly one 1. Let the set of all $n^2 \times n^2$ S-permutation matrices be denoted by $\Sigma_{n^2}$.

As it is proved in \cite{dahl} the cardinality of the set of all S-permutation matrices is equal to
\begin{equation}\label{broy}
\left| \Sigma_{n^2} \right| =  \left( n! \right)^{2n} .
\end{equation}

Two binary   matrices $A=(a_{ij} )\in \mathfrak{B}_{m}$ and $B=( b_{ij} )\in \mathfrak{B}_{m}$ will be called \emph{disjoint} if there are not elements with one and the same indices $a_{ij}$ and $b_{ij}$ such that $a_{ij} =b_{ij} =1$, i.e. if $a_{ij} =1$ then $b_{ij} =0$ and if $b_{ij} =1$ then $a_{ij} =0$, $1\le i,j\le m$.

The concept of S-permutation matrix was introduced by Geir Dahl  \cite{dahl} in relation to the popular Sudoku puzzle. Sudoku is a very popular game. On the other hand, it is well known that Sudoku matrices are special cases of Latin squares in the class of gerechte designs \cite{Bailey}.

Obviously a square $n^2 \times n^2$ matrix $M$ with elements of $[n^2 ] =\{ 1,2,\ldots ,n^2 \}$ is a Sudoku matrix if and only if there are  matrices $A_1 ,A_2 ,\ldots ,A_{n^2} \in\Sigma_{n^2}$, each two of them are disjoint and such that $P$ can be given in the following way:
\begin{equation}\label{disj}
M=1\cdot A_1 +2\cdot A_2 +\cdots +n^2 \cdot A_{n^2}
\end{equation}

Some algorithms for obtaining random Sudoku matrices and their valuation are described in detail in  \cite{yordzhev_random}.

In \cite{Fontana} Roberto Fontana offers an algorithm which randomly gets a family of
$n^2 \times n^2$ mutually disjoint S-permutation matrices, where $n = 2, 3$. In $n = 3$ he
ran the algorithm 1000 times and found 105 different families of nine mutually
disjoint S-permutation matrices. Then using (\ref{disj}) he obtained $9! \cdot 105 = 38\; 102\; 400$ Sudoku matrices.

But it is known~\cite{Felgenhauer} that the total number of $9\times 9$ Sudoku matrices  is
$$9! \cdot 72^2 \cdot 2^7 \cdot 27\; 704\; 267\; 971 = 6\; 670\; 903\; 752\; 021\; 072\; 936\; 960 $$

Thus, in relation with Fontana's algorithm, it looks useful to calculate the probability of two randomly generated S-permutation matrices to be disjoint.
So the question of enumerating all disjoint pairs of  S-permutation matrices naturally arises. This work is devoted to this task.

\emph{Bipartite graph} is the ordered triplet
 $$
g=\langle R_g \cup C_g , E_g \rangle ,
$$
where $R_g$ and $C_g$ are non-empty sets such that $R_g \cap C_g =\emptyset$, the elements of which will be called \emph{vertices}. $E_g \subseteq R_g \times C_g =\{ \langle r,c \rangle \; |\; r\in R_g ,c\in C_g \}$  - the set of \emph{edges}. Multiple edges are not allowed in our considerations.

For more details on graph theory see  \cite{diestel,harary}.

Let $n$ and $k$  be integers such that  $0\le k\le n^2$. Let us denote with $\mathfrak{G}_{n,k}$  the set of all bipartite graphs without multiple edges of the type  $g=\langle  R_g , C_g ,E_g \rangle$, such that  $ |R_g |= |C_g |= n$ and $|E_g |=k$. With $\displaystyle \mathfrak{G}_{n} = \bigcup_{k=0}^{n^2} \mathfrak{G}_{n,k}$  we denote the set of all bipartite graphs without multiple edges of the type  $g=\langle  R_g , C_g ,E_g \rangle $, such that $ |R_g |= |C_g |= n$  and irrespective of the number of the edges.

We will not take into consideration the nature of the vertices of the bipartite graphs, i.e. we will consider that for each of the examined bipartite graphs $g\in \mathfrak{G}_{n,k} $, the vertices are respectively
$$R_g =\{ r_1 ,r_2 ,\ldots ,r_n \}$$
and
$$C_g =\{ c_1 ,c_2 ,\ldots ,c_n \} .$$

Let $g'=\langle R_{g'} \cup C_{g'} , E_{g'} \rangle$ and $g''=\langle R_{g''} \cup C_{g''} , E_{g''} \rangle $ be two bipartite graphs. We will say that $g'$  and $g''$   are \emph{isomorphic} and we will write  $g'\cong g''$, if $|R_{g'}|=|R_{g''}|$, $|C_{g'}|=|C_{g''}|$  and there exist bijections $\rho : R_{g'} \to R_{g''}$ and $ \sigma : C_{g'} \to C_{g''}$, such that  $\langle r, c\rangle \in E_{g'} \Longleftrightarrow \langle \rho (r) ,\sigma (c)\rangle \in E_{g''}$.

Obviously the so entered relation $\cong$  is an equivalence relation and with $\overline{\mathfrak{G}}_{n,k}$  we will denote the factor set (the set of the equivalence classes)

$$\overline{\mathfrak{G}}_{n,k} = {\mathfrak{G}_{n,k}}_{/_{\cong}}$$

Let $g=\langle R_g , C_g ,E_g \rangle \in \mathfrak{G}_{n,k}$  for some natural numbers  $n$ and  $k$. Then with  $\overline{g}$ we will denote the set
$$\overline{g} =\left\{ h\in {\mathfrak{G}}_{n,k} \; |\; h\cong g\right\} \in \overline{\mathfrak{G}}_{n,k}$$
of all isomorphic to $g$  bipartite graphs, and with $|\overline{g}|$  their number. i.e. the cardinality of the set  $\overline{g}$.

Let $g=\langle R_g , C_g ,E_g \rangle \in \mathfrak{G}_{n,k}$  for some natural numbers  $n$  and $k$  and let  $v\in V_g =R_g \cup C_g$.  With  $\gamma (v)$ we will denote the set of all vertices from  $V_g$, adjacent to  $v$, i.e. $u\in \gamma (v)$  if and only if there exists an edge in $E_g$   which joins $u$  and  $v$. If  $v$ is an isolated vertex (i.e. there is no edge incident to  $v$), then by definition  $\gamma (v)=\emptyset$ and  $|\gamma (v)|=0$. Obviously if  $v\in R_g$, then  $\gamma (v)\subseteq C_g$, and if  $v\in C_g$, then  $\gamma (v)\subseteq R_g$.

Obviously
$$\sum_{v\in V_{g}} |\gamma (v)|=2n.$$

Let $n$  and $k$  be positive integers and let  $g\in \mathfrak{G}_{n,k}$. We examine the ordered  $(n+1)$-tuple
$$\langle \psi \rangle  (g)=\langle \psi_0 (g) ,\psi_1 (g),\ldots ,\psi_n (g)\rangle ,$$
where  $\psi_i (g)$, $i=0,1,\ldots ,n$   is equal to the number of the vertices of   $g$, incident to exactly $i$  number of edges. It is easy to see that for each  $g\in \mathfrak{G}_{n,k}$ the equalities $$\displaystyle \sum_{i=1}^n \psi_i (g)=2n,\quad \sum_{i=1}^n i\psi_i (g)=2k .$$ have been executed.

For two bipartite graphs $g$   and  $h$, if  $g\cong h$, then obviously  $$\langle \psi \rangle  (g)=\langle \psi \rangle  (h)=\langle \psi \rangle  (\overline{g})=\langle \psi \rangle  (\overline{h}) .$$

In the set $\mathfrak{B}_n$   we enter equivalence relation  $"\sim "$, such that if $A,B\in \mathfrak{B}_n$  then  $A\sim B$, if  $B$ is obtained from  $A$ after dislocation of some of the rows and/or columns of  $A$.

Let  $g=\langle R_g , C_g , E_g \rangle \in \mathfrak{G}_n$, where $R_g =\{ r_1 ,r_2 ,\ldots ,r_n \}$  and  $C_g =\{ c_1 ,c_2 ,\ldots ,c_n \}$. Then we build the matrix  $A=[a_{ij} ] \in \mathfrak{B}_n$, such that $a_{ij} =1$  if and only if  $\langle r_i ,c_j \rangle \in E_g$.

Inversely, let  $A=[a_{ij} ]\in\mathfrak{B}_n$. We denote the  $i$-th row of $A$  with $r_i$, while the  $j$-th column of $A$  with  $c_j$. Then we build the bipartite graph  $g=\langle R_g , C_g , E_g \rangle \in \mathfrak{G}_n$, where $R_g =\{ r_1 ,r_2 ,\ldots ,r_n \}$, $C_g =\{ c_1 ,c_2 ,\ldots ,c_n \}$  and there exists an edge from the vertex $r_i$  to the vertex $c_j$  if and only if  $a_{ij} =1$.

Thus we showed the following obvious relation between the bipartite graphs and the binary matrices:

\begin{proposition}\label{varphighgh}

There exists one-to-one mapping
$$\varphi : \mathfrak{G}_n \to \mathfrak{B}_n $$
between the elements of $\mathfrak{G}_n$  and  $\mathfrak{B}_n$, such that if $g,h\in \mathfrak{G}_n$, then
$$g\cong h \Longleftrightarrow \varphi (g) \sim \varphi (h) .$$

\end{proposition}

\section{A representation of S-permutation matrices}

If $z_1 \; z_2 \; \ldots \; z_n$ is a permutation of the elements of the set $[n] =\left\{ 1,2,\ldots ,n\right\}$ and we shortly denote $\rho$ this permutation, then in this case we denote by $\rho (i )$ the $i$-th element of this permutation, i.e. $\rho (i) =z_i$, $i=1,2,\ldots ,n$.

Let $\Pi_n$ denotes the set of all $n\times n$ matrices, constructed such that $\pi\in\Pi_n$ if and only if the following three conditions are true:

i) the elements of $\pi$ are ordered pairs of numbers $\langle i,j\rangle$, where $1\le i,j\le n$;

ii) if
$$\left[ \langle a_1 , b_1 \rangle \quad \langle a_2 ,b_2 \rangle \quad \cdots  \quad \langle a_n ,b_n \rangle \right]$$
is the $i$-th row of $\pi$ for any $i\in [n] =\{ 1,2,\ldots ,n\}$, then  $a_1 \; a_2 \; \ldots \; a_n$ in this order is a permutation of the elements of the set $[n]$

iii) if
$$\left[
\begin{array}{c}
\langle a_1 ,b_1 \rangle \\
\langle a_2 ,b_2 \rangle \\
\vdots \\
\langle a_n ,b_n \rangle \\
\end{array}
\right]
$$
is the $j$-th column of $\pi$ for any $j\in [n]$, then  $b_1 ,b_2 ,\ldots , b_n$ in this order is a permutation of the elements of the set $[n]$.

From the definition it follows that every row and every column of any matrix of the set $\Pi_n$ can be identified with permutation of elements of the set $[n]$. Conversely for every $(2n)$-tuple $\langle \langle \rho_1 ,\rho_2 ,\ldots ,\rho_n \rangle ,\langle \sigma_1 ,\sigma_2 ,\ldots , \sigma_n \rangle \rangle$, where $\rho_i = \rho_i (1)\; \rho_i (2) \; \ldots \; \rho_i (n)$, $\sigma_j = \sigma_j (1)\; \sigma_j (2)\; \ldots \; \sigma_j (n)$, $1\le i,j\le n$ are permutations of elements of $[n]$, then the matrix
$$
\pi =
\left[
\begin{array}{cccc}
\langle \rho_1 (1),\sigma_1 (1)\rangle & \langle \rho_1 (2),\sigma_2 (1)\rangle & \cdots & \langle \rho_1 (n),\sigma_n (1)\rangle \\
\langle \rho_2 (1),\sigma_1 (2)\rangle & \langle \rho_2 (2),\sigma_2 (2)\rangle & \cdots & \langle \rho_2 (n),\sigma_n (2)\rangle \\
\vdots & \vdots & \ddots & \vdots \\
\langle \rho_n (1),\sigma_1 (n)\rangle  & \langle \rho_n (2),\sigma_2 (n)\rangle & \cdots & \langle \rho_n (n),\sigma_n (n)\rangle
\end{array}
\right]
$$
is matrix of $\Pi_n$. Hence
\begin{equation}\label{|Pin|}
\left| \Pi_n \right| =\left( n! \right)^{2n}
\end{equation}

\begin{theorem}\label{l2}
Let $n$ be a positive integer, $n\ge 2$. Then there is one to one correspondence between the sets $\Sigma_{n^2}$ and $\Pi_n$.
\end{theorem}

Proof. Let $A\in \Sigma_{n^2}$. Then $A$ is constructed with the help of formula (\ref{matrA}) and for every $i,j\in [n]$ in the block $A_{ij} $ there is only one 1 and let this 1 has coordinates $(a_i ,b_j )$. For every $i,j\in [n]$ we obtain ordered pairs of numbers $\langle a_i ,b_j \rangle$ corresponding to these coordinates. As in every row and every column of $A$ there is only one 1, then the matrix $\left[ \alpha_{ij} \right]_{n\times n}$, where $\alpha_{ij} =\langle a_i ,b_j \rangle $, $1\le i,j\le n$, which is obtained by the ordered pairs of numbers is matrix of $\Pi_n$, i.e. matrix for which the conditions i), ii) and iii) are true.

Conversely, let $\left[ \alpha_{ij} \right]_{n\times n} \in \Pi_n$, where $\alpha_{ij} =\langle a_i ,b_j \rangle $, $i,j \in [n]$, $a_i ,b_j \in [n]$. Then for every $i,j\in [n]$ we construct binary $n\times n$ matrices $A_{ij}$ with only one 1 with coordinates $(a_i ,b_j )$. Then we obtain the matrix of type: (\ref{matrA}). According to the properties i), ii) and iii), it is obvious that the obtained matrix is S-permutation matrix.

\hfill $\Box$

From Theorem \ref{l2} and formula (\ref{|Pin|}) it follows the proof of Proposition \ref{disj} in \cite{dahl}, i.e. formula (\ref{broy}).

\begin{definition}
We say that matrices $\pi ' ,\pi '' \in\Pi_n$, where $\pi ' =\left[ {p'}_{ij} \right]_{n\times n}$, $\pi '' =\left[ {p''}_{ij} \right]_{n\times n}$ are \emph{disjoint}, if ${p'}_{ij} \ne {p''}_{ij}$ for every pair of indices $i,j\in[n]$.
\end{definition}

\begin{proposition}
The number of all  pairs of disjoint matrices of $\Sigma_{n^2}$ is equal to the number of all pairs of disjoint matrices of $\Pi_n$.
\end{proposition}

Proof. It is easy to see that with respect of the described in Theorem \ref{l2} one to one correspondence, every pair of disjoint matrices of $\Sigma_{n^2}$ will correspond to a pair of disjoint matrices of $\Pi_n$ and conversely every pair of disjoint matrices of $\Pi_n$ will correspond to a pair of disjoint matrices of $\Sigma_{n^2}$.

\hfill $\Box$

\begin{definition}
Let $\pi ' ,\pi '' \in\Pi_n$, $\pi ' =\left[ {p'}_{ij} \right]_{n\times n}$, $\pi '' =\left[ {p''}_{ij} \right]_{n\times n}$ and let  the integers $i,j\in[n]$ are such that ${p'}_{ij} = {p''}_{ij}$. In this case we will say that   ${p'}_{ij}$ and ${p''}_{ij}$ are \emph{component-wise equal  elements}.
\end{definition}

Obviously two $\Pi_n$-matrices are disjoint if and only if they do not have component-wise equal elements.

\begin{example}\label{ex1}
\rm We consider the following $\Pi_3$-matrices:
\end{example}

$$
\pi' =\left[ p_{ij}' \right] =
\left[
\begin{array}{ccc}
\langle 1,2\rangle & \langle 3,1\rangle & \langle 2,3\rangle \\
\langle 2,1\rangle & \langle 3,3\rangle & \langle 1,2\rangle \\
\langle 3,3\rangle & \langle 1,2\rangle & \langle 2,1\rangle
\end{array}
\right]
$$

$$
\pi'' =\left[ p_{ij}'' \right] =
\left[
\begin{array}{ccc}
\langle 1,3\rangle & \langle 3,2\rangle & \langle 2,1\rangle \\
\langle 3,1\rangle & \langle 1,1\rangle & \langle 2,2\rangle \\
\langle 3,2\rangle & \langle 1,3\rangle & \langle 2,3\rangle
\end{array}
\right]
$$

$$
\pi''' =\left[ p_{ij}''' \right] =
\left[
\begin{array}{ccc}
\langle 1,2\rangle & \langle 3,3\rangle & \langle 2,1\rangle \\
\langle 2,1\rangle & \langle 3,2\rangle & \langle 1,2\rangle \\
\langle 3,3\rangle & \langle 1,1\rangle & \langle 2,3\rangle
\end{array}
\right]
$$

Matrices $\pi'$ and $\pi''$ are disjoint, because they do not have component-wise equal elements.

Matrices $\pi''$ and $\pi'''$ are not disjoint, because they have two component-wise equal elements: $p_{13}' =p_{13}''' =\langle 2,1\rangle$ and $p_{33}'' =p_{33}''' =\langle 2,3\rangle$.

Matrices $\pi'$ and $\pi'''$ are not disjoint, because they have four component-wise equal elements: $p_{11}' =p_{11}''' =\langle 1,2\rangle$, $p_{21}' =p_{21}''' =\langle 2,1\rangle$, $p_{23}' =p_{23}''' =\langle 1,2\rangle$ and $p_{31}' =p_{31}''' =\langle 3,3\rangle$.

\section{A  formula  for counting  disjoint pairs  of $n^2 \times n^2$  S-permutation matrices}

\begin{lemma}\label{l3dskmat}
Let  $\pi\in\Pi_n$. Then the number $q(n,k)$ of all matrices $\pi'  \in\Pi_n$,  having at least  $k$, $k=0,1,\ldots ,n^2$ component-wise equal elements to the matrix  $\pi$ is equal to
\begin{equation}\label{fl3gfgfg}
q(n,k)= \sum_{\overline{g}\in \overline{\mathfrak{G}}_{n,k} } |\overline{g} | \left( \prod_{i=0}^{n-2} \left[ \left( n-i\right) ! \right]^{\psi_i (\overline{g})} \right)
\end{equation}
\end{lemma}

Proof. Let  $\pi =\left[ p_{ij} \right]_{n\times n} ,\pi' =\left[ p'_{ij} \right]_{n\times n} \in \Pi_n$  and let $\pi$  and $\pi'$  have exactly $k$ component-wise equal elements. Then we uniquely  obtain the binary $n\times n$  matrix  $A=\left[ a_{ij} \right]_{n\times n}$, such that $a_{ij} =1$  if and only if $p_{ij} =p'_{ij}$, $i,j\in [n]$. Let graph $g\in \mathfrak{G}_{n,k}$  be such that  $g=\varphi^{-1} (A)$, where $\varphi :\mathfrak{B}_n \to \mathfrak{G}_n$  is the one-to-one mapping which gave us the grounds to formulate Proposition \ref{varphighgh}. This graph $g$ identically corresponds to the ordered pair of matrices  $\langle \pi ,\pi' \rangle \in \Pi\times\Pi$.

Inversely, let $g=\langle R_g , C_g ,E_g \rangle \in \mathfrak{G}_{n,k}$, $V_g =R_g \cup C_g$, $A=[a_{ij} ]_{n\times n} =\varphi (g)$ and let $\pi =\left[ p_{ij} \right]_{n\times n}$  be an arbitrary matrix from  $\Pi_n$. We search for the number $h(\pi ,A)$ of all matrices  $\pi' =[p'_{ij} ]_{n\times n} \in \Pi_n$, such that  $p'_{ij} =p_{ij}$, if  $a_{ij}=1$. (It is assumed that there exist $s,t\in [n]$  such that $a_{st} =0$  and  $p'_{st} =p_{st}$.) Let the  $i$-th row, $i=1,2,\ldots ,n$  of $\pi$  correspond to the permutation  $\rho_i$ of the elements of $[n]$  and let the  $i$-th row of the matrix  $A$ correspond to the vertex $r_i \in R_g$  of graph  . Then there exist $(n-|\gamma (r_i )|)!$  permutations $\rho'_i$  of the elements of  $ [n]$, such that if   $a_{it} =1$, then $\rho_i (t) =\rho'_i (t)$, $t\in [n]$. Likewise we also prove the respective statement for the columns of   $\pi$. Therefore
$$\displaystyle h(\pi ,A) =\prod_{v\in V_g} (n-|\gamma (v)|)! .$$

From everything said so far it follows that for each $\pi\in\Pi_n$  there exist
$$q(n,k)=\sum_{g\in\mathfrak{G}_{n,k}} \left( \prod_{v\in V_g} (n-|\gamma (v)|)! \right)$$
matrices from  $\Pi_n$, which have at least $k$ elements that are component-wise equal to the respective elements of $\pi$.

But obviously, if $h\in\mathfrak{G}_{n,k}$  is such a bipartite graph that  $h\cong g$, then
$$\displaystyle \left( \prod_{v\in V_g} (n-|\gamma (g)|)! \right) = \left( \prod_{u\in V_h} (n-|\gamma (h)|)! \right) = \left( \prod_{i=0}^{n} \left[ \left( n-i\right) ! \right]^{\psi_i (\overline{g})} \right),$$	
whence follows the equality
$$q(n,k)= \sum_{\overline{g}\in \overline{\mathfrak{G}}_{n,k} } |\overline{g} | \left( \prod_{v\in R_g \cup C_g} ( n-|\gamma (v)|)! \right) = \sum_{\overline{g}\in \overline{\mathfrak{G}}_{n,k} } |\overline{g} | \left( \prod_{i=0}^{n} \left[ \left( n-i\right) ! \right]^{\psi_i (\overline{g})} \right)$$

And since $(n-n)!=0!=1$  and  $[n-(n-1)]!=1!=1$, then we finally obtain formula (\ref{fl3gfgfg}).

\hfill $\Box$

\begin{theorem}\label{mainTh}
Let  $A\in \Sigma_{n^2}$.Then the number $\xi_n$ of all matrices $B\in \Sigma_{n^2}$ which are disjoint with $A$ is equal to
\begin{equation}\label{gl6_main}
\xi_{n} = (n!)^{2n} +\sum_{k=1}^{n^2} \left( -1\right)^k \left( \sum_{\overline{g}\in \overline{\mathfrak{G}}_{n,k} } |\overline{g} | \left( \prod_{i=0}^{n-2} \left[ \left( n-i\right) ! \right]^{\psi_i (\overline{g})} \right) \right)
\end{equation}
\end{theorem}

Proof. Let $n\ge 2$  be an integer. Then applying theorem \ref{l2}, lemma \ref{l3dskmat} and the principle of inclusion and exclusion we obtain that the number $\xi_n$  of all matrices $B\in \Sigma_{n^2}$  which are disjoint with  $A$ is equal to
$$\xi_n =   |\Pi_n | +\sum_{k=1}^{n^2} (-1)^k q(n,k),$$
where the function $q(n,k)$  is calculated with the help of formula (\ref{fl3gfgfg}), while $|\Pi_n |$  with the help of formula (\ref{|Pin|}). Thus we obtain the proof to formula (\ref{gl6_main}).

\hfill $\Box$

\begin{corollary}\label{th2-gl6}
The cardinality $\eta_{n}$ of the set of all  disjoint non-ordered  pairs of $n^2 \times n^2$ S-permutation matrices is equal to
\begin{equation}\label{nonordereddisjointpair_gl6}
\eta_{n} =\frac{(n!)^{2n}}{2} \xi_n
\end{equation}
where $\xi_n$ is described using formula \ref{gl6_main}.
\end{corollary}

Proof. follows directly from formula (\ref{broy}) and having in mind that the ''disjoint'' relation is symmetric and antireflexive.

\hfill $\Box$

\begin{corollary}\label{th3_gl6}
The probability $p(n)$ of two randomly generated $n^2 \times n^2$ S-permutation matrices to be disjoint is equal to
\begin{equation}\label{probbility_gl6}
p(n) = \frac{\displaystyle \xi_n}{\displaystyle  \left( n! \right)^{2n} -1}  ,
\end{equation}
where $\xi_n$ is described using formula \ref{gl6_main}.
\end{corollary}

Proof. Applying Corollary \ref{th2-gl6} and formula (\ref{broy}), we obtain:

$$p(n)= \frac{\displaystyle \eta_{n}}{\displaystyle {\left| \Sigma_{n^2} \right| \choose 2}} = \frac{\displaystyle \frac{(n!)^{2n}}{2} \xi_n}{\displaystyle \frac{\left( n! \right)^{2n} \left( \left( n! \right)^{2n} -1\right) }{2}} =$$
$$ = \frac{\displaystyle \xi_n}{\displaystyle  \left( n! \right)^{2n} -1} .$$

\hfill $\Box$

\section{Calculation of the Number of the Disjoint Pairs of S-permutation matrices when $n=2$ and $n=3$}\label{ISRN1}
\subsection{Consider $n=2$}
When   $n=2$, $\overline{\mathfrak{G}}_{2}$ is composed of seven equivalence classes also including the graph without edges ($k=0$), which does not participate in our calculations. When  $k=1,2,3,4$, we have depicted one representative $g_1$, $g_2$, $g_3$, $g_4$, $g_5$  and $g_6$  from each equivalence class respectively on figures  \ref{n2k1}, \ref{n2k2}, \ref{n2k3} and \ref{n2k4}.

\unitlength=0.7mm \linethickness{0.5pt}

\begin{figure}
\begin{picture}(50,30)

\put(0,0.01){\framebox(50,30)}
\put(3,27){\makebox(0,0){$g_1$}}

\multiput(13,16)(24,0){2}{\oval(13,20)}
\put(14,3){\makebox(0,0){$R_{g_1}$}}
\put(38,3){\makebox(0,0){$C_{g_1}$}}

\put(13,20){\circle*{2}}
\put(37,20){\circle{2}}
\put(15,20){\line(1,0){20}}

\put(13,10){\circle*{2}}
\put(37,10){\circle{2}}

\end{picture}
\caption{$n=2$, $k=1$}\label{n2k1}
\end{figure}

\begin{figure}
\begin{picture}(170,30)

\multiput(0,0)(60,0){3}{
\put(0,0.01){\framebox(50,30)}

\multiput(13,16)(24,0){2}{\oval(13,20)}

\put(13,20){\circle*{2}}
\put(37,20){\circle{2}}

\put(13,10){\circle*{2}}
\put(37,10){\circle{2}}
}

\put(15,20){\line(1,0){20}}
\put(15,10){\line(1,0){20}}

\put(75,20){\line(1,0){20}}
\put(75,10){\line(2,1){20}}

\put(135,20){\line(1,0){20}}
\put(135,20){\line(2,-1){20}}

\put(3,27){\makebox(0,0){$g_2$}}
\put(63,27){\makebox(0,0){$g_3$}}
\put(123,27){\makebox(0,0){$g_4$}}

\put(14,3){\makebox(0,0){$R_{g_2}$}}
\put(38,3){\makebox(0,0){$C_{g_2}$}}

\put(74,3){\makebox(0,0){$R_{g_3}$}}
\put(98,3){\makebox(0,0){$C_{g_3}$}}

\put(134,3){\makebox(0,0){$R_{g_4}$}}
\put(158,3){\makebox(0,0){$C_{g_4}$}}

\end{picture}
\caption{$n=2$, $k=2$}\label{n2k2}
\end{figure}

\begin{figure}
\begin{picture}(50,30)

\put(0,0.01){\framebox(50,30)}
\put(3,27){\makebox(0,0){$g_5$}}

\multiput(13,16)(24,0){2}{\oval(13,20)}
\put(14,3){\makebox(0,0){$R_{g_5}$}}
\put(38,3){\makebox(0,0){$C_{g_5}$}}

\put(13,20){\circle*{2}}
\put(37,20){\circle{2}}

\put(15,20){\line(2,-1){20}}

\put(13,10){\circle*{2}}
\put(37,10){\circle{2}}
\put(15,10){\line(1,0){20}}
\put(15,10){\line(2,1){20}}

\end{picture}
\caption{$n=2$, $k=3$}\label{n2k3}
\end{figure}

\begin{figure}
\begin{picture}(50,30)

\put(0,0.01){\framebox(50,30)}
\put(3,27){\makebox(0,0){$g_6$}}

\multiput(13,16)(24,0){2}{\oval(13.5,20)}
\put(14,3){\makebox(0,0){$R_{g_6}$}}
\put(38,3){\makebox(0,0){$C_{g_6}$}}

\put(13,20){\circle*{2}}
\put(37,20){\circle{2}}
\put(15,20){\line(1,0){20}}
\put(15,20){\line(2,-1){20}}

\put(13,10){\circle*{2}}
\put(37,10){\circle{2}}
\put(15,10){\line(1,0){20}}
\put(15,10){\line(2,1){20}}

\end{picture}
\caption{$n=2$, $k=4$}\label{n2k4}
\end{figure}

It is not difficult to notice that\\
$\displaystyle |\overline{g}_1 |=4,\quad  \langle \psi \rangle  (\overline{g}_1 ) =\langle 2,2,0\rangle ,\quad k=1$\\
$\displaystyle |\overline{g}_2 |=2,\quad \langle \psi \rangle  (\overline{g}_2 ) =\langle 0,4,0\rangle ,\quad k=2$\\
$\displaystyle |\overline{g}_3 |=|\overline{g}_4 |=2,\quad \langle \psi \rangle  (\overline{g}_3 ) =\langle \psi \rangle  (\overline{g}_4 ) =\langle 1,2,1\rangle  ,\quad k=2$\\
$\displaystyle |\overline{g}_5 |=4,\quad \langle \psi \rangle  (\overline{g}_5 )=\langle 0,2,2\rangle ,\quad k=3$\\
$\displaystyle |\overline{g}_6 |=1,\quad \langle \psi \rangle  (\overline{g}_6 )=\langle 0,0,4\rangle ,\quad k=4$

Let  $A\in \Sigma_{4}$. Then the number $\xi_2$ of all matrices  $B\in \Sigma_{4}$ which are disjoint with  $A$ is equal to

$$\xi_{2} = (2!)^{4} +\sum_{k=1}^{4} \left( -1\right)^k \left( \sum_{\overline{g}\in \overline{\mathfrak{G}}_{n,k} } |\overline{g} | \cdot 2^{\psi_0 (\overline{g})}  \right) =$$
$$=16-4\cdot 2^2 +\left( 2\cdot 2^0 +2\cdot 2^1 +2\cdot 2^1 \right) -4\cdot 2^0 +1\cdot 2^0 =7$$

$$\eta_{2} =\frac{2^{4}}{2} \xi_2 =56$$

$$p(2) = \frac{\displaystyle \xi_2}{\displaystyle   2^{4} -1} =\frac{7}{15}$$

\subsection{Consider $n=3$}

On figures from \ref{n3k1} to  \ref{n3k9} one representative from each equivalence class of the factor sets $\mathfrak{\overline{G}}_{3,k}$, $k=1,2,\ldots ,9$ has been depicted, and we have numbered these graphs from 7 to 41.

When $n=3$  formula (\ref{gl6_main}) has the form of
$$\xi_3 = 6^6 + \sum_{i=7}^{41} |\overline{g}_i | 6^{\psi_0 (\overline{g}_i )} 2^{\psi_1 (\overline{g}_i )} (-1)^{\kappa(\overline{g}_i )} \  ,$$
where with $\kappa(\overline{g}_i )$ we have denoted the number of the edges of the graphs from the equivalence class  $\overline{g}_i$.

\begin{figure}
\begin{picture}(50,40)
\put(0,0){\framebox(50,40)}
\put(3,37){\makebox(0,0){$g_7$}}
\multiput(13,21)(24,0){2}{\oval(14,30)}
\put(14,3){\makebox(0,0){$R_{g_7}$}}
\put(38,3){\makebox(0,0){$C_{g_7}$}}
\put(13,30){\circle*{2}}
\put(37,30){\circle{2}}
\put(15,30){\line(1,0){20}}
\put(13,20){\circle*{2}}
\put(37,20){\circle{2}}
\put(13,10){\circle*{2}}
\put(37,10){\circle{2}}
\end{picture}
\caption{$n=3$, $k=1$}\label{n3k1}
\end{figure}

\begin{figure}
\begin{picture}(170,40)

\multiput(0,0)(60,0){3}{
\put(0,0){\framebox(50,40)}

\multiput(13,21)(24,0){2}{\oval(14,30)}

\put(13,30){\circle*{2}}
\put(37,30){\circle{2}}

\put(13,20){\circle*{2}}
\put(37,20){\circle{2}}

\put(13,10){\circle*{2}}
\put(37,10){\circle{2}}
}

\put(3,37){\makebox(0,0){$g_8$}}
\put(63,37){\makebox(0,0){$g_9$}}
\put(123,37){\makebox(0,0){$g_{10}$}}

\put(14,3){\makebox(0,0){$R_{g_8}$}}
\put(38,3){\makebox(0,0){$C_{g_8}$}}
\put(74,3){\makebox(0,0){$R_{g_9}$}}
\put(98,3){\makebox(0,0){$C_{g_9}$}}
\put(134,3){\makebox(0,0){$R_{g_{10}}$}}
\put(158,3){\makebox(0,0){$C_{g_{10}}$}}

\put(15,30){\line(1,0){20}}
\put(15,20){\line(1,0){20}}

\put(75,30){\line(1,0){20}}
\put(75,30){\line(2,-1){20}}

\put(135,30){\line(1,0){20}}
\put(135,20){\line(2,1){20}}

\end{picture}
\caption{$n=3$, $k=2$}\label{n3k2}
\end{figure}

\begin{figure}
\begin{picture}(170,90)

\multiput(0,0)(60,0){3}{
\multiput(0,0)(0,50){2}{

\put(0,0){\framebox(50,40)}

\multiput(13,21)(24,0){2}{\oval(14,30)}

\put(13,30){\circle*{2}}
\put(37,30){\circle{2}}

\put(13,20){\circle*{2}}
\put(37,20){\circle{2}}

\put(13,10){\circle*{2}}
\put(37,10){\circle{2}}
}
}

\put(3,87){\makebox(0,0){$g_{11}$}}
\put(63,87){\makebox(0,0){$g_{12}$}}
\put(123,87){\makebox(0,0){$g_{13}$}}

\put(14,53){\makebox(0,0){$R_{g_{11}}$}}
\put(38,53){\makebox(0,0){$C_{g_{11}}$}}

\put(74,53){\makebox(0,0){$R_{g_{12}}$}}
\put(98,53){\makebox(0,0){$C_{g_{12}}$}}

\put(134,53){\makebox(0,0){$R_{g_{13}}$}}
\put(158,53){\makebox(0,0){$C_{g_{13}}$}}

\put(3,37){\makebox(0,0){$g_{14}$}}
\put(63,37){\makebox(0,0){$g_{15}$}}
\put(123,37){\makebox(0,0){$g_{16}$}}

\put(14,3){\makebox(0,0){$R_{g_{14}}$}}
\put(38,3){\makebox(0,0){$C_{g_{14}}$}}

\put(74,3){\makebox(0,0){$R_{g_{15}}$}}
\put(98,3){\makebox(0,0){$C_{g_{15}}$}}

\put(134,3){\makebox(0,0){$R_{g_{16}}$}}
\put(158,3){\makebox(0,0){$C_{g_{16}}$}}

\put(15,80){\line(1,0){20}}
\put(15,70){\line(1,0){20}}
\put(15,60){\line(1,0){20}}

\put(75,80){\line(1,0){20}}
\put(75,80){\line(2,-1){20}}
\put(75,60){\line(1,0){20}}

\put(135,80){\line(1,0){20}}
\put(135,70){\line(2,1){20}}
\put(135,60){\line(1,0){20}}

\put(15,30){\line(1,0){20}}
\put(15,30){\line(2,-1){20}}
\put(15,20){\line(1,0){20}}

\put(75,20){\line(1,0){20}}
\put(75,20){\line(2,1){20}}
\put(75,20){\line(2,-1){20}}

\put(135,30){\line(2,-1){20}}

\put(135,20){\line(1,0){20}}

\put(135,10){\line(2,1){20}}

\end{picture}
\caption{$n=3$, $k=3$}\label{n3k3}
\end{figure}

\begin{figure}
\begin{picture}(170,140)

\multiput(0,0)(60,0){2}{
\multiput(0,0)(0,50){3}{

\put(0,0){\framebox(50,40)}

\multiput(13,21)(24,0){2}{\oval(14,30)}

\put(13,30){\circle*{2}}
\put(37,30){\circle{2}}

\put(13,20){\circle*{2}}
\put(37,20){\circle{2}}

\put(13,10){\circle*{2}}
\put(37,10){\circle{2}}
}
}
\put(120,100){\framebox(50,40)}

\multiput(133,121)(24,0){2}{\oval(14,30)}

\put(133,130){\circle*{2}}
\put(157,130){\circle{2}}

\put(133,120){\circle*{2}}
\put(157,120){\circle{2}}

\put(133,110){\circle*{2}}
\put(157,110){\circle{2}}

\put(3,137){\makebox(0,0){$g_{17}$}}
\put(63,137){\makebox(0,0){$g_{18}$}}
\put(123,137){\makebox(0,0){$g_{19}$}}

\put(14,103){\makebox(0,0){$R_{g_{17}}$}}
\put(38,103){\makebox(0,0){$C_{g_{17}}$}}

\put(74,103){\makebox(0,0){$R_{g_{18}}$}}
\put(98,103){\makebox(0,0){$C_{g_{18}}$}}

\put(134,103){\makebox(0,0){$R_{g_{19}}$}}
\put(158,103){\makebox(0,0){$C_{g_{19}}$}}

\put(3,87){\makebox(0,0){$g_{20}$}}
\put(63,87){\makebox(0,0){$g_{21}$}}

\put(14,53){\makebox(0,0){$R_{g_{20}}$}}
\put(38,53){\makebox(0,0){$C_{g_{20}}$}}

\put(74,53){\makebox(0,0){$R_{g_{21}}$}}
\put(98,53){\makebox(0,0){$C_{g_{21}}$}}

\put(3,37){\makebox(0,0){$g_{22}$}}
\put(63,37){\makebox(0,0){$g_{23}$}}

\put(14,3){\makebox(0,0){$R_{g_{22}}$}}
\put(38,3){\makebox(0,0){$C_{g_{22}}$}}

\put(74,3){\makebox(0,0){$R_{g_{23}}$}}
\put(98,3){\makebox(0,0){$C_{g_{23}}$}}

\put(15,130){\line(1,0){20}}
\put(15,130){\line(2,-1){20}}

\put(15,120){\line(1,0){20}}
\put(15,120){\line(2,1){20}}

\put(75,130){\line(1,0){20}}
\put(75,130){\line(2,-1){20}}

\put(75,120){\line(2,-1){20}}

\put(75,110){\line(1,0){20}}

\put(135,130){\line(1,0){20}}

\put(135,130){\line(2,-1){20}}

\put(135,120){\line(1,0){20}}

\put(135,110){\line(1,0){20}}

\put(15,80){\line(1,0){20}}

\put(15,70){\line(2,1){20}}
\put(15,70){\line(2,-1){20}}

\put(15,60){\line(1,0){20}}

\put(75,80){\line(1,0){20}}

\put(75,80){\line(2,-1){20}}

\put(75,60){\line(1,0){20}}

\put(75,60){\line(2,1){20}}

\put(15,30){\line(1,0){20}}

\put(15,20){\line(1,0){20}}
\put(15,20){\line(2,1){20}}
\put(15,20){\line(2,-1){20}}

\put(75,30){\line(1,0){20}}

\put(75,30){\line(2,-1){20}}

\put(75,20){\line(1,0){20}}

\put(75,10){\line(2,1){20}}

\end{picture}
\caption{$n=3$, $k=4$}\label{n3k4}
\end{figure}

\begin{figure}
\begin{picture}(170,140)

\multiput(0,0)(60,0){2}{
\multiput(0,0)(0,50){3}{

\put(0,0){\framebox(50,40)}

\multiput(13,21)(24,0){2}{\oval(14,30)}

\put(13,30){\circle*{2}}
\put(37,30){\circle{2}}

\put(13,20){\circle*{2}}
\put(37,20){\circle{2}}

\put(13,10){\circle*{2}}
\put(37,10){\circle{2}}
}
}
\put(120,100){\framebox(50,40)}

\multiput(133,121)(24,0){2}{\oval(14,30)}

\put(133,130){\circle*{2}}
\put(157,130){\circle{2}}

\put(133,120){\circle*{2}}
\put(157,120){\circle{2}}

\put(133,110){\circle*{2}}
\put(157,110){\circle{2}}

\put(3,137){\makebox(0,0){$g_{24}$}}
\put(63,137){\makebox(0,0){$g_{25}$}}
\put(123,137){\makebox(0,0){$g_{26}$}}

\put(14,103){\makebox(0,0){$R_{g_{24}}$}}
\put(38,103){\makebox(0,0){$C_{g_{24}}$}}

\put(74,103){\makebox(0,0){$R_{g_{25}}$}}
\put(98,103){\makebox(0,0){$C_{g_{25}}$}}

\put(134,103){\makebox(0,0){$R_{g_{26}}$}}
\put(158,103){\makebox(0,0){$C_{g_{26}}$}}

\put(3,87){\makebox(0,0){$g_{27}$}}
\put(63,87){\makebox(0,0){$g_{28}$}}

\put(14,53){\makebox(0,0){$R_{g_{27}}$}}
\put(38,53){\makebox(0,0){$C_{g_{27}}$}}

\put(74,53){\makebox(0,0){$R_{g_{28}}$}}
\put(98,53){\makebox(0,0){$C_{g_{28}}$}}

\put(3,37){\makebox(0,0){$g_{29}$}}
\put(63,37){\makebox(0,0){$g_{30}$}}

\put(14,3){\makebox(0,0){$R_{g_{29}}$}}
\put(38,3){\makebox(0,0){$C_{g_{29}}$}}

\put(74,3){\makebox(0,0){$R_{g_{30}}$}}
\put(98,3){\makebox(0,0){$C_{g_{30}}$}}

\put(15,130){\line(1,-1){20}}

\put(15,120){\line(2,-1){20}}

\put(15,110){\line(1,0){20}}
\put(15,110){\line(1,1){20}}
\put(15,110){\line(2,1){20}}

\put(75,130){\line(1,-1){20}}

\put(75,120){\line(1,0){20}}
\put(75,120){\line(2,1){20}}

\put(75,110){\line(1,1){20}}
\put(75,110){\line(2,1){20}}

\put(135,130){\line(1,-1){20}}

\put(135,120){\line(2,1){20}}
\put(135,120){\line(2,-1){20}}

\put(135,110){\line(1,1){20}}
\put(135,110){\line(2,1){20}}

\put(15,80){\line(1,-1){20}}
\put(15,80){\line(2,-1){20}}

\put(15,70){\line(1,0){20}}

\put(15,60){\line(1,1){20}}
\put(15,60){\line(2,1){20}}

\put(75,80){\line(1,-1){20}}

\put(75,70){\line(1,0){20}}
\put(75,70){\line(2,1){20}}
\put(75,70){\line(2,-1){20}}

\put(75,60){\line(1,1){20}}

\put(15,30){\line(1,-1){20}}
\put(15,30){\line(2,-1){20}}

\put(15,10){\line(1,0){20}}
\put(15,10){\line(1,1){20}}
\put(15,10){\line(2,1){20}}

\put(75,30){\line(1,-1){20}}

\put(75,20){\line(2,1){20}}
\put(75,20){\line(2,-1){20}}

\put(75,10){\line(1,0){20}}
\put(75,10){\line(1,1){20}}

\end{picture}
\caption{$n=3$, $k=5$}\label{n3k5}
\end{figure}

\begin{figure}
\begin{picture}(170,90)

\multiput(0,0)(60,0){3}{
\multiput(0,0)(0,50){2}{

\put(0,0){\framebox(50,40)}

\multiput(13,21)(24,0){2}{\oval(14,30)}

\put(13,30){\circle*{2}}
\put(37,30){\circle{2}}

\put(13,20){\circle*{2}}
\put(37,20){\circle{2}}

\put(13,10){\circle*{2}}
\put(37,10){\circle{2}}
}
}

\put(3,87){\makebox(0,0){$g_{31}$}}
\put(63,87){\makebox(0,0){$g_{32}$}}
\put(123,87){\makebox(0,0){$g_{33}$}}

\put(14,53){\makebox(0,0){$R_{g_{31}}$}}
\put(38,53){\makebox(0,0){$C_{g_{31}}$}}

\put(74,53){\makebox(0,0){$R_{g_{32}}$}}
\put(98,53){\makebox(0,0){$C_{g_{32}}$}}

\put(134,53){\makebox(0,0){$R_{g_{33}}$}}
\put(158,53){\makebox(0,0){$C_{g_{33}}$}}

\put(3,37){\makebox(0,0){$g_{34}$}}
\put(63,37){\makebox(0,0){$g_{35}$}}
\put(123,37){\makebox(0,0){$g_{36}$}}

\put(14,3){\makebox(0,0){$R_{g_{34}}$}}
\put(38,3){\makebox(0,0){$C_{g_{34}}$}}

\put(74,3){\makebox(0,0){$R_{g_{35}}$}}
\put(98,3){\makebox(0,0){$C_{g_{35}}$}}

\put(134,3){\makebox(0,0){$R_{g_{36}}$}}
\put(158,3){\makebox(0,0){$C_{g_{36}}$}}

\put(15,80){\line(1,-1){20}}
\put(15,80){\line(2,-1){20}}

\put(15,70){\line(2,1){20}}
\put(15,70){\line(2,-1){20}}

\put(15,60){\line(1,1){20}}
\put(15,60){\line(2,1){20}}

\put(75,80){\line(1,-1){20}}

\put(75,70){\line(1,0){20}}
\put(75,70){\line(2,1){20}}
\put(75,70){\line(2,-1){20}}

\put(75,60){\line(1,1){20}}
\put(75,60){\line(2,1){20}}

\put(135,80){\line(1,-1){20}}
\put(135,80){\line(2,-1){20}}

\put(135,70){\line(1,0){20}}
\put(135,70){\line(2,-1){20}}

\put(135,60){\line(1,1){20}}
\put(135,60){\line(2,1){20}}

\put(15,30){\line(1,-1){20}}

\put(15,20){\line(2,1){20}}
\put(15,20){\line(2,-1){20}}

\put(15,10){\line(1,0){20}}
\put(15,10){\line(1,1){20}}
\put(15,10){\line(2,1){20}}

\put(75,30){\line(1,0){20}}
\put(75,30){\line(1,-1){20}}
\put(75,30){\line(2,-1){20}}

\put(75,10){\line(1,0){20}}
\put(75,10){\line(1,1){20}}
\put(75,10){\line(2,1){20}}

\put(135,30){\line(1,0){20}}
\put(135,30){\line(1,-1){20}}

\put(135,20){\line(2,1){20}}
\put(135,20){\line(2,-1){20}}

\put(135,10){\line(1,0){20}}
\put(135,10){\line(1,1){20}}

\end{picture}
\caption{$n=3$, $k=6$}\label{n3k6}
\end{figure}

\begin{figure}
\begin{picture}(170,40)

\multiput(0,0)(60,0){3}{
\put(0,0){\framebox(50,40)}

\multiput(13,21)(24,0){2}{\oval(14,30)}

\put(13,30){\circle*{2}}
\put(37,30){\circle{2}}

\put(13,20){\circle*{2}}
\put(37,20){\circle{2}}

\put(13,10){\circle*{2}}
\put(37,10){\circle{2}}
}

\put(3,37){\makebox(0,0){$g_{37}$}}
\put(63,37){\makebox(0,0){$g_{38}$}}
\put(123,37){\makebox(0,0){$g_{39}$}}

\put(14,3){\makebox(0,0){$R_{g_{37}}$}}
\put(38,3){\makebox(0,0){$C_{g_{37}}$}}
\put(74,3){\makebox(0,0){$R_{g_{38}}$}}
\put(98,3){\makebox(0,0){$C_{g_{38}}$}}
\put(134,3){\makebox(0,0){$R_{g_{39}}$}}
\put(158,3){\makebox(0,0){$C_{g_{39}}$}}

\put(15,30){\line(2,-1){20}}
\put(15,30){\line(1,-1){20}}
\put(15,20){\line(2,1){20}}
\put(15,20){\line(2,-1){20}}
\put(15,10){\line(1,0){20}}
\put(15,10){\line(1,1){20}}
\put(15,10){\line(2,1){20}}

\put(75,30){\line(1,-1){20}}
\put(75,20){\line(2,1){20}}
\put(75,20){\line(1,0){20}}
\put(75,20){\line(2,-1){20}}
\put(75,10){\line(1,0){20}}
\put(75,10){\line(1,1){20}}
\put(75,10){\line(2,1){20}}

\put(135,30){\line(2,-1){20}}
\put(135,30){\line(1,-1){20}}
\put(135,20){\line(1,0){20}}
\put(135,20){\line(2,-1){20}}
\put(135,10){\line(1,0){20}}
\put(135,10){\line(1,1){20}}
\put(135,10){\line(2,1){20}}

\end{picture}
\caption{$n=3$, $k=7$}\label{n3k7}
\end{figure}

\begin{figure}
\begin{picture}(50,40)

\put(0,0){\framebox(50,40)}
\put(3,37){\makebox(0,0){$g_{40}$}}

\multiput(13,21)(24,0){2}{\oval(14,30)}
\put(14,3){\makebox(0,0){$R_{g_{40}}$}}
\put(38,3){\makebox(0,0){$C_{g_{40}}$}}

\put(13,30){\circle*{2}}
\put(37,30){\circle{2}}
\put(15,30){\line(1,-1){20}}
\put(15,30){\line(2,-1){20}}

\put(13,20){\circle*{2}}
\put(37,20){\circle{2}}
\put(15,20){\line(1,0){20}}
\put(15,20){\line(2,1){20}}
\put(15,20){\line(2,-1){20}}

\put(13,10){\circle*{2}}
\put(37,10){\circle{2}}
\put(15,10){\line(1,0){20}}
\put(15,10){\line(1,1){20}}
\put(15,10){\line(2,1){20}}

\end{picture}
\caption{$n=3$, $k=8$}\label{n3k8}
\end{figure}

\begin{figure}
\begin{picture}(50,40)

\put(0,0){\framebox(50,40)}
\put(3,37){\makebox(0,0){$g_{41}$}}

\multiput(13,21)(24,0){2}{\oval(14,30)}
\put(14,3){\makebox(0,0){$R_{g_{41}}$}}
\put(38,3){\makebox(0,0){$C_{g_{41}}$}}

\put(13,30){\circle*{2}}
\put(37,30){\circle{2}}
\put(15,30){\line(1,0){20}}
\put(15,30){\line(1,-1){20}}
\put(15,30){\line(2,-1){20}}

\put(13,20){\circle*{2}}
\put(37,20){\circle{2}}
\put(15,20){\line(1,0){20}}
\put(15,20){\line(2,1){20}}
\put(15,20){\line(2,-1){20}}

\put(13,10){\circle*{2}}
\put(37,10){\circle{2}}
\put(15,10){\line(1,0){20}}
\put(15,10){\line(1,1){20}}
\put(15,10){\line(2,1){20}}

\end{picture}
\caption{$n=3$, $k=9$}\label{n3k9}
\end{figure}

Below we  enumerate the examined numerical characteristics of the respective equivalence classes of bipartite graphs. Their calculation is trivial.\\
$\displaystyle |\overline{g}_7 |={3\choose 1}^2 =9,\quad \langle \psi \rangle  (\overline{g}_7 ) =\langle 4,2,0,0\rangle ,\quad k=1$\\
$\displaystyle |\overline{g}_8 |={3\choose 1}^2 |\overline{g}_2 |=18,\quad  \langle \psi \rangle  (\overline{g}_8 ) =\langle 2,4,0,0\rangle ,\quad k=2$\\
$\displaystyle |\overline{g}_9 |=|\overline{g}_{10} |={3\choose 1}{3\choose 2} =9,\quad \langle \psi \rangle  (\overline{g}_9 )= \langle \psi \rangle  (\overline{g}_{10} ) =\langle 3,2,1,0\rangle ,\quad k=2$
$\displaystyle |\overline{g}_{11} |= |\mathcal{S}_3 |=3! =6,\quad \langle \psi \rangle  (\overline{g}_{11} ) =\langle 0,6,0,0\rangle ,\quad k=3$\\
$\displaystyle |\overline{g}_{12} |=|\overline{g}_{13} |={3\choose 1}^2 |\overline{g}_3 |=18,\quad \langle \psi \rangle  (\overline{g}_{12} ) =\langle \psi \rangle  (\overline{g}_{13} ) =\langle 1,4,1,0\rangle ,\quad k=3$\\
$\displaystyle |\overline{g}_{14} |={3\choose 1}^2 |\overline{g}_5 |=36,\quad \langle \psi \rangle  (\overline{g}_{14} ) =\langle 2,2,2,0\rangle ,\quad k=3$\\
$\displaystyle |\overline{g}_{15} |=|\overline{g}_{16} |={3\choose 1}=3,\quad \langle \psi \rangle  (\overline{g}_{15} ) =\langle \psi \rangle  (\overline{g}_{16} ) =\langle 2,3,0,1\rangle ,\quad k=3$\\
$\displaystyle |\overline{g}_{17} |={3\choose 1}^2 =9,\quad \langle \psi \rangle  (\overline{g}_{17} ) =\langle 2,0,4,0\rangle ,\quad k=4$\\
$\displaystyle |\overline{g}_{18} |={3\choose 1}{3\choose2} =9,\quad \langle \psi \rangle  (\overline{g}_{18} ) =\langle 0,4,2,0\rangle ,\quad k=4$\\
$\displaystyle |\overline{g}_{19} |={3\choose 1}^2 |\overline{g}_5 |=36,\quad \langle \psi \rangle  (g_{19} ) =\langle 0,4,2,0\rangle,\quad k=4$\\
$\displaystyle |\overline{g}_{20} |=|\overline{g}_{21} |={3\choose 1}^2 |\overline{g}_2 |=18,\quad \langle \psi \rangle  (\overline{g}_{20} ) =\langle \psi \rangle  (\overline{g}_{21} ) =\langle 1,2,3,0\rangle ,\quad k=4$\\
$\displaystyle |\overline{g}_{22} |=|\overline{g}_{23} |={3\choose 1}^2 {2\choose 1} =18 ,\quad \langle \psi \rangle  (\overline{g}_{22} ) =\langle \psi \rangle  (\overline{g}_{23} ) =\langle 1,3,1,1\rangle ,\ k=4$\\
$\displaystyle |\overline{g}_{24} |={3\choose 1}^2 =9,\quad \langle \psi \rangle  (\overline{g}_{24} ) =\langle 0,4,0,2\rangle ,\quad k=5$\\
$\displaystyle |\overline{g}_{25} |={3\choose 1}^2 =9,\quad \langle \psi \rangle  (g_{25} ) =\langle 0,2,4,0\rangle ,\quad k=5$\\
$\displaystyle |\overline{g}_{26} |={3\choose 1}^2 {2\choose 1}^2 =36,\quad \langle \psi \rangle  (\overline{g}_{26} ) =\langle 0,2,4,0\rangle ,\quad k=5$\\
$\displaystyle |\overline{g}_{27} |=|\overline{g}_{28} |={3\choose 1}^2 |\overline{g}_{2} |=18,\quad \langle \psi \rangle  (g_{27} ) =\langle \psi \rangle  (g_{28} ) =\langle 0,3,2,1\rangle ,\quad k=5$\\
$\displaystyle |\overline{g}_{29} |=|\overline{g}_{30} |=3! {3\choose 1} =18,\quad \langle \psi \rangle  (\overline{g}_{29} ) =\langle \psi \rangle  (\overline{g}_{30} ) =\langle 1,1,3,1\rangle ,\quad k=5$\\
$\displaystyle |\overline{g}_{31} |=|\overline{g}_{11} |=6,\quad \langle \psi \rangle  (\overline{g}_{31} ) =\langle 0,0,6,0\rangle ,\quad k=6$\\
$\displaystyle |\overline{g}_{32} |=|\overline{g}_{33} |=|\overline{g}_{12} |=18,\quad \langle \psi \rangle  (\overline{g}_{32} ) =\langle \psi \rangle  (\overline{g}_{33} ) =\langle 0,1,4,1\rangle ,\quad k=6$\\
$\displaystyle |\overline{g}_{34} |=|\overline{g}_{14} |=36,\quad \langle \psi \rangle  (\overline{g}_{34} ) =\langle 0,2,2,2\rangle ,\quad k=6$\\
$\displaystyle |\overline{g}_{35} |=|\overline{g}_{36} |=|\overline{g}_{15} |=3,\quad \langle \psi \rangle  (g_{35} ) =\langle \psi \rangle  (g_{36} ) =\langle 1,0,3,2\rangle ,\quad k=6$\\
$\displaystyle |\overline{g}_{37} |=|\overline{g}_8 |=18,\quad \displaystyle \langle \psi \rangle  (g_{37} ) =\langle 0,0,4,2\rangle ,\quad k=7$\\
$\displaystyle |\overline{g}_{38} |=|\overline{g}_{39} |=|\overline{g}_9 |=9,\quad \langle \psi \rangle  (g_{38} ) =\langle \psi \rangle  (g_{39} ) =\langle 0,1,2,3\rangle ,\quad k=7$\\
$\displaystyle |\overline{g}_{40} |={3\choose 1}^2 =9,\quad \langle \psi \rangle  (g_{40} ) =\langle 0,0,2,4\rangle ,\quad k=8$\\
$\displaystyle |\overline{g}_{41} |=1,\quad  \langle \psi \rangle  (g_{41} ) =\langle 0,0,0,6\rangle ,\quad k=9$\\

In order to calculate  $\xi_3$  we used the next computer programme written in programming language C++. 

\begin{verbatim}
int main() {
	double xi_3 =pow(6,6);
	int g[][4] = {
		9,4,2,1,
		18,2,4,2,
		9,3,2,2,
		9,3,2,2,
		6,0,6,3,
		18,1,4,3,
		18,1,4,3,
		36,2,2,3,
		3,2,3,3,
		3,2,3,3,
		9,2,0,4,
		9,0,4,4,
		36,0,4,4,
		18,1,2,4,
		18,1,2,4,
		18,1,3,4,
		18,1,3,4,
		9,0,4,5,
		9,0,2,5,
		36,0,2,5,
		18,0,3,5,
		18,0,3,5,
		18,1,1,5,
		18,1,1,5,
		6,0,0,6,
		18,0,1,6,
		18,0,1,6,
		36,0,2,6,
		3,1,0,6,
		3,1,0,6,
		18,0,0,7,
		9,0,1,7,
		9,0,1,7,
		9,0,0,8,
		1,0,0,9
	};
	for (int i=0; i<35; i++)
		xi_3 += g[i][0] * pow(6,g[i][1]) * pow(2,g[i][2]) * (g[i][3]%2 ? -1 : 1);
	double p3;
	long eta_3;
	eta_3 = (long) (pow(6,6) /2) * xi_3;
	p3 =xi_3 / (pow(6,6)-1);
	cout<<"xi_3 = "<<xi_3<<endl;
	cout<<"eta_3 = "<<eta_3<<endl;
	cout<<"p(3) = "<<p3<<endl;
	return 0;
}
\end{verbatim}

After its work we obtain that
$$\xi_3 =17\; 972 $$

Then according to formula (\ref{nonordereddisjointpair_gl6}), the number of all disjoint nonordered pairs of $9 \times 9$  S-permutation matrices is equal to
$$\eta_{3} =\frac{(3!)^{6}}{2} \xi_3 =  419\; 250\; 816.$$

The probability  $p(3)$ for two randomly obtained $9 \times 9$  S-permutation matrices to be disjoint, according to formula (\ref{probbility_gl6}) is equal to
$$p(3) = \frac{\displaystyle \xi_3}{\displaystyle  \left( 3! \right)^{6} -1}  =       0.385211$$

\bibliographystyle{plain}
\bibliography{S-permMatr}

\end{document}